\documentclass{article}

\begin{document}

\vspace{0.5cm}

\centerline{\LARGE \bf Note on the Algebra of Screening Currents }
\centerline{\LARGE \bf for the Quantum Deformed $W$-Algebra}

\vspace{0.5cm}

\centerline{Liu Zhao\hspace{1cm} Bo-Yu Hou}

\vspace{0.2cm}

\centerline{\it Institute of Modern Physics, Northwest University, Xian 710069, China}
\centerline{\it e-mail: lzhao@nwu.edu.cn,\hspace{0.3cm}byhou@nwu.edu.cn}

\vspace{0.6cm}

\begin{abstract}
With slight modifications in the zero modes contributions, the positive
and negative screening currents for the quantum deformed $W$-algebra
${\cal W}_{q,p}(g)$ can be put together to form a single algebra which
can be regarded as an elliptic deformation of the universal enveloping
algebra of $\widehat{g}$, where $g$ is any classical simply-laced Lie algebra.
\end{abstract}

\vspace{1cm}


Recently, various deformations of the classical and quantum Virasoro and
$W$-algebras have received considerable interests. N.Reshetikhin and E.Frenkel \cite{a}
first introduced the Poisson algebras ${\cal W}_q(g)$, which are $q$-deformation
of classical $W$-algebras. Later on, J.Shiraishi, H.Kubo, H.Awata and S.Odake \cite{b}
obtained a quantum version of the algebra ${\cal W}_q(sl_2)$, which is a
noncommutative algebra depending on two parameters $p$ and $q$.
H.Awata, H.Kubo, S.Odake and J.Shiraishi \cite{awata}
and B.Feigin, E.Frenkel \cite{c}
independently extended this result to general case, i.e.
quantum deformed $W$-algebras
${\cal W}_{q,p}(g)$, where $g$ is any classical semisimple Lie algebra.
All these algebras were obtained together with their respective bosonic Fock
space representations. Similar considerations with respect
to the Yangian deformation were also carried out and have led to $\hbar$
deformed Virasoro algebra \cite{d} and quantum $(\xi,\hbar)$-deformed
$W$-algebras \cite{e}.

In their work \cite{c}, B.Feigin and E.Frenkel also
obtained the screening currents for the
algebra ${\cal W}_{q,p}(g)$ and found the elliptic relations between them.
They noticed that the relations for the positive (resp. negative) screening currents
form an elliptic deformation for the loop algebra $\hat{\sf n}_+$
(resp. $\hat{\sf n}_-$) of the nilpotent subalgebra ${\sf n}_+$
(resp. ${\sf n}_-$) of $g$. However, they did not consider whether these two nilpotent
elliptic algebras can be put together to form a unified elliptic algebra.

In this note we shall show that it is possible to combine the above nilpotent
algebras into a single unified elliptic algebra if we introduce some new generating
currents denoted by $H^\pm_i(z)$ and slightly modify the zero mode contributions
in the bosonic representation of the screening currents (the modified
``screening currents''\footnote{Strictly speaking, the modified currents are
no longer screening currents of the deformed $W$-algebra and hence we here use
the quotation marks.} are denoted by $E_i(z)$ and
$F_i(z)$ respectively whilst the original ones by $S^\pm_i(z)$).
Unlike the unmodified screening currents, the modified currents
$E_i(z)$, $F_i(z)$ and the newly introduced currents $H^\pm_i(z)$ do not commute
with the $W$-algebra generating currents even up to total differences.

Let us start from a brief description of the results of B.Feigin and E.Frenkel
\cite{c} which are necessary for our discussion. First we consider the simple case of
$g=sl_N$. By definition, the algebra ${\cal W}_{q,p}(sl_N)$ is generated by the
Fourier coefficients of the currents $T_1(z), \cdots, T_{N-1}(z)$, which, in the
free field realization, obey the quantum deformed Miura transformation

\begin{eqnarray*}
& &D_{p^{-1}}^N - T_1(z) D_{p^{-1}}^{N-1} + T_2(z) D_{p^{-1}}^{N-2} -
 \cdots + (-1)^{N-1} T_{N-1}(z) D_{p^{-1}} + (-1)^N\\
& &~~~~~~~= :(D_{p^{-1}}- \Lambda_1(z)) (D_{p^{-1}}- \Lambda_2(zp))
\cdots (D_{p^{-1}}- \Lambda_N(zp^{N-1})):,
\end{eqnarray*}

\noindent where $\Lambda_i, i=1, \cdots, N$ are generating series of a Heisenberg
algebra, $[ D_a f ](x) =f(xa)$. The screening currents $S^\pm_i(z)$
are solutions of the difference equations

\begin{eqnarray*}
D_{q} S^+_i(z) = p^{-1}:\Lambda_{i+1}(zp^{i/2})\Lambda_i(zp^{i/2})^{-1}S^+_i(z):,\\
D_{p/q} S^-_i(z) = p^{-1}:\Lambda_{i+1}(zp^{i/2})\Lambda_i(zp^{i/2})^{-1}S^-_i(z):.
\end{eqnarray*}

\noindent From the above formulas one can show that the screening currents
$S^+_i(z)$ and $S^-_i(z)$ commute with the $W$-algebra generating currents
$T_i(z)$ up to a total difference. The case of other simply-laced
$g$ of rank $N-1$ can be understood in a similar fashion with appropriate
modification in the form of Miura transformation.

Let $p,q$ be two generic parameters such that $|p/q|<1$. We introduce a
third parameter $\beta$ using the relation $p=q^{1-\beta}$.
For general simply-laced Lie algebra $g$ with Cartan matrix $(A_{ij})$, we
introduce the Heisenberg algebra \cite{c} ${\cal H}_{q,p}(g)$ with generators
$a_i[n], ~i=1, \cdots, N-1; n\in Z$ and $Q_i, i=1, \cdots, N-1$ and relations

\begin{eqnarray*}
& &{}[ a_i[n], a_j[m] ] =\frac{1}{n}\frac{(1-q^n)(p^{A_{ij}n/2} - p^{-A_{ij}n/2})
(1-(p/q)^n)}{1-p^n} \delta_{n,-m},\\
& &{}[a_i[n], Q_j] = A_{ij}\beta \delta_{n,0}.
\end{eqnarray*}

\noindent Let

\begin{eqnarray*}
& & s_i^+[m] = \frac{a_i[m]}{q^{-m}-1},~~~m\neq 0,~~~~s^+_i[0]=a_i[0],\\
& & s_i^-[m] = -\frac{a_i[m]}{(q/p)^{m}-1},~~~m\neq 0,~~~~s^-_i[0]=a_i[0]/\beta.
\end{eqnarray*}

\noindent Then the screening currents $S^+_i(z)$ and $S^-_i(z)$
can be realized in the Fock space of ${\cal H}_{q,p}(g)$ as \cite{c}

\begin{eqnarray}
& & S^+_i(z)=\mbox{e}^{Q_i}z^{s^+_i[0]}
:\exp\left(\sum_{m\neq0}s^+_i[m]z^{-m}\right):,\\
& & S^-_i(z)=\mbox{e}^{-Q_i/\beta}z^{-s^-_i[0]}
:\exp\left(-\sum_{m\neq0}s^-_i[m]z^{-m}\right):.
\end{eqnarray}

Using these bosonic expressions, I.Feigin and E.Frenkel \cite{c}
arrived at the following
elliptic relations for the screening currents $S^+_i(z)$ and $S^-_i(z)$,

\begin{eqnarray}
& &S^+_i(z) S^+_j(w) = (-1)^{A_{ij}-1} \left(\frac{w}{z}\right)^{A_{ij}-A_{ij}\beta-1}
\frac{\theta_q\left(\frac{w}{z}p^{A_{ij}/2}\right)}
{\theta_q\left(\frac{z}{w}p^{A_{ij}/2}\right)} S^+_j(w) S^+_i(z), \label{7}\\
& &S^-_i(z) S^-_j(w) = (-1)^{A_{ij}-1} \left(\frac{w}{z}\right)^{A_{ij}-A_{ij}/\beta-1}
\frac{\theta_{p/q}\left(\frac{w}{z}p^{A_{ij}/2}\right)}
{\theta_{p/q}\left(\frac{z}{w}p^{A_{ij}/2}\right)} S^-_j(w) S^-_i(z), \label{8}
\end{eqnarray}

\noindent where as usual,

\begin{eqnarray*}
\theta_a(x) =(x|a)_\infty (ax^{-1}| a)_\infty (a|a)_\infty,~~~~
(x|a)\infty \equiv \prod_{n=0}^\infty (1-xa^n),
\end{eqnarray*}

\noindent and the equations (\ref{7},\ref{8}) are to be understood in the sense of
analytical continuation. Notice that the function $\theta_a(x)$ is an
elliptic function with multiplicative periods $a$ and $\mbox{e}^{2\pi i}$,

\begin{eqnarray*}
\theta_a(ax) = -x^{-1} \theta_a(x),~~~~~\theta_a(x\mbox{e}^{2\pi i}) =\theta_a(x).
\end{eqnarray*}

\noindent In \cite{c}, I.Feigin and E.Frenkel also obtained the cross relations
between $S^+_i(z)$ and $S^-_i(z)$ as normal ordered relations,

\begin{eqnarray}
& &S^+_i(z) S^-_i(w) = \frac{1}{(z-wq)(z-wp^{-1}q)}:S^+_i(z) S^-_i(w):, \label{10}\\
& &S^+_i(z) S^-_j(w) = (z-wp^{-1/2}q) :S^+_i(z)S^-_j(w):,~~~A_{ij}=-1,\\
& &S^+_i(z) S^-_j(w) = :S^+_i(z) S^-_j(w):,~~~A_{ij}=0.
\end{eqnarray}

\noindent The reversed relation can also be obtained straightforwardly,

\begin{eqnarray}
& &S^-_i(w) S^+_i(z) = \frac{1}{(w-zq^{-1})(w-zpq^{-1})}:S^+_i(z) S^-_i(w):,\\
& &S^-_j(w) S^+_i(z) = (w-zp^{1/2}q^{-1}) :S^+_i(z)S^-_j(w):,~~~A_{ij}=-1,\\
& &S^-_j(w) S^+_i(z) = :S^+_i(z) S^-_j(w):,~~~A_{ij}=0. \label{15}
\end{eqnarray}

The primary motivation of this work was to combine the algebra of positive and negative
screening currents $S^\pm_i(z)$ into a single unified algebra. For this the
relation between $S^+_i(z)$ and $S^-_j(w)$ have to be closed in the sense of
commutator algebra. However, as can be easily shown from eqs.(\ref{10})-(\ref{15}), this
is impossible if we use the original form of the screening currents, because
in the right hand side of the commutation relation (here $\simeq$ implies ``equals
up to regular terms'')

\begin{eqnarray*}
{}[ S^+_i(z), S^-_j(w) ] \simeq \delta_{ij}
\left(\frac{1}{(z-wq)(z-wp^{-1}q)}-\frac{1}{(w-zq^{-1})(w-zpq^{-1})}\right)
:S^+_i(z)S^-_j(w):,
\end{eqnarray*}

\noindent the operator $:S^+_i(z)S^-_j(w):$ is bi-localized and the rational
expression $\frac{1}{(z-wq)(z-wp^{-1}q)}-\frac{1}{(w-zq^{-1})(w-zpq^{-1})}$
cannot be rewritten as the sum of $\delta$-functions. This imply that the commutator
$[ S^+_i(z), S^-_j(w) ]$ does not close over the space of operators depending
only on one spectral parameter.

The way around this difficulty is to modify the screening currents slightly so that the
commutator of the modified ``screening currents'' indeed close over the space of operators
depending on a single parameter.

The modified ``screening currents'' read

\begin{eqnarray}
& & E_i(z)=\mbox{e}^{Q_i}(z(p/q)^{1/2})^{P_i}
:\exp\left(\sum_{m\neq0}s^+_i[m]z^{-m}\right):, \label{17}\\
& & F_i(z)=\mbox{e}^{-Q_i} (zq^{1/2})^{-P_i}
:\exp\left(-\sum_{m\neq0}s^-_i[m]z^{-m}\right):, \label{18}
\end{eqnarray}

\noindent where $P_i=s^-_i[0]$ which satisfy the relation

\begin{eqnarray*}
{}[P_i, Q_j] = A_{ij}.
\end{eqnarray*}

After the above modification, we can easily show that the relations for $E_i(z)$
and $F_i(z)$ read

\begin{eqnarray}
& &E_i(z) E_j(w) = (-1)^{A_{ij}-1} \left(\frac{w}{z}\right)^{-1}
\frac{\theta_q\left(\frac{w}{z}p^{A_{ij}/2}\right)}
{\theta_q\left(\frac{z}{w}p^{A_{ij}/2}\right)} E_j(w) E_i(z),\label{19}\\
& &F_i(z) F_j(w) = (-1)^{A_{ij}-1} \left(\frac{w}{z}\right)^{-1}
\frac{\theta_{p/q}\left(\frac{w}{z}p^{A_{ij}/2}\right)}
{\theta_{p/q}\left(\frac{z}{w}p^{A_{ij}/2}\right)} F_j(w) F_i(z). \label{20}
\end{eqnarray}

\noindent Moreover, equations (\ref{10})-(\ref{15}) will be turned into

\begin{eqnarray}
& &E_i(z) F_i(w) = \frac{1}{\left(z(p/q)^{1/2}\right)^2 (1-\frac{wq}{z})
(1-\frac{wp^{-1}q}{z})}:E_i(z) F_i(w):,\\
& &E_i(z) E_j(w) = \left(z(p/q)^{1/2}\right)(1-\frac{w}{z}p^{-1/2}q)
:E_i(z)F_j(w):,~~~A_{ij}=-1,\\
& &E_i(z) F_j(w) = :E_i(z) F_j(w):,~~~A_{ij}=0,\\
& &F_i(w) E_i(z) = \frac{1}{\left(wq^{1/2}\right)^2 (1-\frac{z}{wq})
(1-\frac{z}{wp^{-1}q})}:E_i(z) F_i(w):,\\
& &F_j(w) E_i(z) = \left(wq^{1/2}\right) (1-\frac{z}{w}p^{1/2}q^{-1})
:E_i(z)F_j(w):,~~~A_{ij}=-1,\\
& &F_j(w) E_i(z) = :E_i(z) F_j(w):,~~~A_{ij}=0,
\end{eqnarray}

\noindent following from which we have the commutation relation

\begin{eqnarray}
{}[ E_i(z), F_j(w) ] \simeq \frac{\delta_{ij}}{(p-1)zw}
\left[ \delta(\frac{w}{zq}) H^+_i(zq^{-1/2}) - \delta(\frac{w}{z(p/q)})
H^-_i(w(p/q)^{-1/2}) \right],  \label{21}
\end{eqnarray}

\noindent where we have introduced the new generating currents $H^\pm_i(z)$

\begin{eqnarray}
& & H^+_i(z) = :E_i(zq^{1/2}) F_i(zq^{-1/2}):, \label{22}\\
& & H^-_i(z) = :E_i(z(p/q)^{-1/2}) F_i(z(p/q)^{1/2}):,   \label{23}
\end{eqnarray}

\noindent which play the role of Cartan subalgebra generators, and the
$\delta$-function is defined as

\begin{eqnarray*}
\delta(z) = \sum_{n\in Z} z^n,~~~~~f(z)\delta(z/w)=f(w)\delta(z/w).
\end{eqnarray*}

\noindent From eqs.(\ref{19})-(\ref{21}) we can further obtain
the following commutation relations
which should all be understood in the sense of analytical continuations,

\begin{eqnarray}
H^\pm_i(z)H^\pm_j(w)&=&\left(\frac{w}{z}\right)^{-2}
\frac{\theta_q\left(\frac{w}{z}p^{A_{ij}/2}\right)
\theta_{\tilde{q}}\left(\frac{w}{z}p^{A_{ij}/2}\right)}
{\theta_q\left(\frac{z}{w}p^{A_{ij}/2}\right)
\theta_{\tilde{q}}\left(\frac{z}{w}p^{A_{ij}/2}\right)}
H^\pm_j(w)H^\pm_i(z), \label{24}\\
H^+_i(z)H^-_j(w)&=&\left(\frac{w}{z}\right)^{-2}
\frac{\theta_q\left(\frac{w}{z}p^{(A_{ij}}p^{-1/2}\right)
\theta_{\tilde{q}}\left(\frac{w}{z}p^{(A_{ij}}p^{1/2}\right)}
{\theta_q\left(\frac{z}{w}p^{(A_{ij}}p^{1/2}\right)
\theta_{\tilde{q}}\left(\frac{z}{w}p^{(A_{ij}}p^{-1/2}\right)}
H^-_j(w)H^+_i(z), \\
H^+_i(z)E_j(w)&=& (-1)^{A_{ij}-1} \left(\frac{w}{zq^{1/2}}\right)^{-1}
\frac{\theta_q\left(\frac{w}{z}p^{A_{ij}/2}q^{-1/2}\right)}
{\theta_q\left(\frac{z}{w}p^{A_{ij}/2}q^{1/2}\right)}
E_j(w)H^+_i(z),\\
H^-_i(z)E_j(w)&=& (-1)^{A_{ij}-1} \left(\frac{w}{z(p/q)^{-1/2}}\right)^{-1}
\frac{\theta_q\left(\frac{w}{z}p^{A_{ij}/2}(p/q)^{1/2}\right)}
{\theta_q\left(\frac{z}{w}p^{A_{ij}/2}(p/q)^{-1/2}\right)}
E_j(w)H^-_i(z),\\
H^+_i(z)F_j(w)&=& (-1)^{A_{ij}-1} \left(\frac{w}{zq^{-1/2}}\right)^{-1}
\frac{\theta_{p/q}\left(\frac{w}{z}p^{A_{ij}/2}q^{1/2}\right)}
{\theta_{p/q}\left(\frac{z}{w}p^{A_{ij}/2}q^{-1/2}\right)}
F_j(w)H^+_i(z),\\
H^-_i(z)F_j(w)&=& (-1)^{A_{ij}-1} \left(\frac{w}{z(p/q)^{1/2}}\right)^{-1}
\frac{\theta_{p/q}\left(\frac{w}{z}p^{A_{ij}/2}(p/q)^{-1/2}\right)}
{\theta_{p/q}\left(\frac{z}{w}p^{A_{ij}/2}(p/q)^{1/2}\right)}
F_j(w)H^-_i(z). \label{29}
\end{eqnarray}

For $g=sl_2$, eqs.(\ref{19}),(\ref{20}),(\ref{21}),(\ref{24})-(\ref{29})
already define an elliptic algebra which can be regarded
as the $c=1$ case of the following more general elliptic deformation of the
universal enveloping algebra of the Kac-Moody algebra $\widehat{sl}_2$
\footnote{We would like to thank H.Awata for bringing our attention
to \cite{awatab}, where the $c=1$ form of the elliptic algebra in the
case of $g=sl_2$ was obtained independently.},

\begin{eqnarray}
H^\pm(z)H^\pm(w)&=&\left(\frac{w}{z}\right)^{-2}
\frac{\theta_q\left(\frac{w}{z}p\right)
\theta_{\tilde{q}}\left(\frac{w}{z}p\right)}
{\theta_q\left(\frac{z}{w}p\right)
\theta_{\tilde{q}}\left(\frac{z}{w}p\right)}
H^\pm(w)H^\pm(z), \\
H^+(z)H^-(w)&=&\left(\frac{w}{z}\right)^{-2}
\frac{\theta_q\left(\frac{w}{z}p^{(2-c)/2}\right)
\theta_{\tilde{q}}\left(\frac{w}{z}p^{(2+c)/2}\right)}
{\theta_q\left(\frac{z}{w}p^{(2+c)/2}\right)
\theta_{\tilde{q}}\left(\frac{z}{w}p^{(2-c)/2}\right)}
H^-(w)H^+(z), \\
H^+(z)E(w)&=& - \left(\frac{wq^{-c/2}}{z}\right)^{-1}
\frac{\theta_q\left(\frac{w}{z}p q^{-c/2}\right)}
{\theta_q\left(\frac{z}{w}p q^{c/2}\right)}
E(w)H^+(z),\\
H^-(z)E(w)&=& - \left(\frac{w\tilde{q}^{c/2}}{z}\right)^{-1}
\frac{\theta_q\left(\frac{w}{z}p \tilde{q}^{c/2}\right)}
{\theta_q\left(\frac{z}{w}p \tilde{q}^{-c/2}\right)}
E(w)H^-(z),\\
H^+(z)F(w)&=& - \left(\frac{wq^{c/2}}{z}\right)^{-1}
\frac{\theta_{\tilde{q}}\left(\frac{w}{z}p q^{c/2}\right)}
{\theta_{\tilde{q}}\left(\frac{z}{w}p q^{-c/2}\right)}
F(w)H^+(z),\\
H^-(z)F(w)&=& - \left(\frac{w\tilde{q}^{-c/2}}{z}\right)^{-1}
\frac{\theta_{\tilde{q}}\left(\frac{w}{z}p \tilde{q}^{-c/2}\right)}
{\theta_{\tilde{q}}\left(\frac{z}{w}p \tilde{q}^{c/2}\right)}
F(w)H^-(z),\\
E(z)E(w)&=& - \left(\frac{w}{z}\right)^{-1}
\frac{\theta_q\left(\frac{w}{z}p \right)}
{\theta_q\left(\frac{z}{w}p \right)}
E(w)E(z), \label{36}\\
F(z)F(w)&=& - \left(\frac{w}{z}\right)^{-1}
\frac{\theta_{\tilde{q}}\left(\frac{w}{z}p \right)}
{\theta_{\tilde{q}}\left(\frac{z}{w}p \right)}
F(w)F(z),\label{37}\\
{}[ E(z), F(w) ] &=&\frac{1}{(p-1)zw} \left[
\delta\left(\frac{z}{wq^c}\right)H^+(zq^{-c/2})
-\delta\left(\frac{w}{z\tilde{q}^c}\right)H^-(w\tilde{q}^{-c/2}) \right],
\end{eqnarray}

\noindent where $c$ is the central charge and $\tilde{q}$ and $q$ are connected
by the relation

\begin{eqnarray*}
q\tilde{q}=p^c.
\end{eqnarray*}

\noindent For general $g$, eqs.(\ref{19}),(\ref{20}),(\ref{21}),(\ref{24})-(\ref{29})
 do not yet form a closed algebra
because the cubic Serre-like relations for $E_i(z)$ (resp. $F_i(z)$) have not been
supplemented. Such Serre-like relations can be explicitly obtained using the
results (\ref{36}) and (\ref{37}).
The final closed algebra has the following generating relations,

\begin{eqnarray}
H^\pm_i(z)H^\pm_j(w)&=&\left(\frac{w}{z}\right)^{-2}
\frac{\theta_q\left(\frac{w}{z}p^{A_{ij}/2}\right)
\theta_{\tilde{q}}\left(\frac{w}{z}p^{A_{ij}/2}\right)}
{\theta_q\left(\frac{z}{w}p^{A_{ij}/2}\right)
\theta_{\tilde{q}}\left(\frac{z}{w}p^{A_{ij}/2}\right)}
H^\pm_j(w)H^\pm_i(z), \label{39}\\
H^+_i(z)H^-_j(w)&=&\left(\frac{w}{z}\right)^{-2}
\frac{\theta_q\left(\frac{w}{z}p^{(A_{ij}-c)/2}\right)
\theta_{\tilde{q}}\left(\frac{w}{z}p^{(A_{ij}+c)/2}\right)}
{\theta_q\left(\frac{z}{w}p^{(A_{ij}+c)/2}\right)
\theta_{\tilde{q}}\left(\frac{z}{w}p^{(A_{ij}-c)/2}\right)}
H^-_j(w)H^+_i(z), \\
H^+_i(z)E_j(w)&=& (-1)^{A_{ij}-1} \left(\frac{wq^{-c/2}}{z}\right)^{-1}
\frac{\theta_q\left(\frac{w}{z}p^{A_{ij}/2} q^{-c/2}\right)}
{\theta_q\left(\frac{z}{w}p^{A_{ij}/2} q^{c/2}\right)}
E_j(w)H^+_i(z),\\
H^-_i(z)E_j(w)&=& (-1)^{A_{ij}-1} \left(\frac{w\tilde{q}^{c/2}}{z}\right)^{-1}
\frac{\theta_q\left(\frac{w}{z}p^{A_{ij}/2} \tilde{q}^{c/2}\right)}
{\theta_q\left(\frac{z}{w}p^{A_{ij}/2} \tilde{q}^{-c/2}\right)}
E_j(w)H^-_i(z),\\
H^+_i(z)F_j(w)&=& (-1)^{A_{ij}-1} \left(\frac{wq^{c/2}}{z}\right)^{-1}
\frac{\theta_{\tilde{q}}\left(\frac{w}{z}p^{A_{ij}/2} q^{c/2}\right)}
{\theta_{\tilde{q}}\left(\frac{z}{w}p^{A_{ij}/2} q^{-c/2}\right)}
F_j(w)H^+_i(z),\\
H^-_i(z)F_j(w)&=& (-1)^{A_{ij}-1} \left(\frac{w\tilde{q}^{-c/2}}{z}\right)^{-1}
\frac{\theta_{\tilde{q}}\left(\frac{w}{z}p^{A_{ij}/2} \tilde{q}^{-c/2}\right)}
{\theta_{\tilde{q}}\left(\frac{z}{w}p^{A_{ij}/2} \tilde{q}^{c/2}\right)}
F_j(w)H^-_i(z),\\
E_i(z)E_j(w)&=& (-1)^{A_{ij}-1} \left(\frac{w}{z}\right)^{-1}
\frac{\theta_q\left(\frac{w}{z}p^{A_{ij}/2} \right)}
{\theta_q\left(\frac{z}{w}p^{A_{ij}/2} \right)}
E_j(w)E_i(z),\\
F_i(z)F_j(w)&=& (-1)^{A_{ij}-1} \left(\frac{w}{z}\right)^{-1}
\frac{\theta_{\tilde{q}}\left(\frac{w}{z}p^{A_{ij}/2} \right)}
{\theta_{\tilde{q}}\left(\frac{z}{w}p^{A_{ij}/2} \right)}
F_j(w)F_i(z),\\
{}[ E_i(z), F_j(w) ] &=&\frac{\delta_{ij}}{(p-1)zw} \left[
\delta\left(\frac{z}{wq^c}\right)H^+(zq^{-c/2})
-\delta\left(\frac{w}{z\tilde{q}^c}\right)H^-(w\tilde{q}^{-c/2}) \right],\\
E_i(z_1)E_i(z_2)E_j(w) \!\!\!\!&-&\!\!\!\! f_{ij}(z_1/w,z_2/w) E_i(z_1)E_j(w)E_i(z_2) +
E_j(w)E_i(z_1)E_i(z_2)\nonumber\\
& & + (\mbox{replacement } ~z_1 \leftrightarrow z_2) =0,
\hspace{0.5cm} A_{ij}=-1,\\
F_i(z_1)F_i(z_2)F_j(w) \!\!\!\!&-&\!\!\!\! g_{ij}(z_1/w,z_2/w) F_i(z_1)F_j(w)F_i(z_2) +
F_j(w)F_i(z_1)F_i(z_2) \nonumber\\
& & + (\mbox{replacement } ~z_1 \leftrightarrow z_2) =0,
\hspace{0.5cm} A_{ij}=-1, \label{51}
\end{eqnarray}

\noindent where

\begin{eqnarray*} 
& &f_{ij}(z_1/w,z_2/w) =
\frac{\left(\psi^{(q)}_{ii}\left(\frac{z_2}{z_1}\right)+1 \right)
\left(\psi^{(q)}_{ij}\left(\frac{w}{z_1}\right)
\psi^{(q)}_{ij}\left(\frac{w}{z_2}\right)+1 \right)}
{\psi^{(q)}_{ij}\left(\frac{w}{z_2}\right)
+\psi^{(q)}_{ii}\left(\frac{z_2}{z_1}\right)
\psi^{(q)}_{ij}\left(\frac{w}{z_1}\right)},\\
& &g_{ij}(z_1/w,z_2/w)
= \frac{\left(\psi^{(\tilde{q})}_{ii}\left(\frac{z_2}{z_1}\right)+1 \right)
\left(\psi^{(\tilde{q})}_{ij}\left(\frac{w}{z_1}\right)
\psi^{(\tilde{q})}_{ij}\left(\frac{w}{z_2}\right)+1 \right)}
{\psi^{(\tilde{q})}_{ij}\left(\frac{w}{z_2}\right)
+\psi^{(\tilde{q})}_{ii}\left(\frac{z_2}{z_1}\right)
\psi^{(\tilde{q})}_{ij}\left(\frac{w}{z_1}\right)},\\
\end{eqnarray*}

\noindent in which

\begin{eqnarray*}
& &\psi^{(q)}_{ij}(x) =  (-1)^{A_{ij}-1} \left(x\right)^{-1}
\frac{\theta_q\left(xp^{A_{ij}/2} \right)}
{\theta_q\left(x^{-1}p^{A_{ij}/2} \right)},\\
& &\psi^{(\tilde{q})}_{ij}(x) =  (-1)^{A_{ij}-1} \left(x\right)^{-1}
\frac{\theta_{\tilde{q}}\left(xp^{A_{ij}/2} \right)}
{\theta_{\tilde{q}}\left(x^{-1}p^{A_{ij}/2} \right)}
\end{eqnarray*}

\noindent are structure functions appeared in the commutation relations
between $E_i(z), E_j(w)$ and $F_i(z), F_j(w)$, and they admit the
factorization property

\begin{eqnarray*}
\psi^{(q)}_{ij}(x) = \frac{\phi^{(q)}_{ij}(x)}{\phi^{(q)}_{ij}(x^{-1})},~~~~
\psi^{(\tilde{q})}_{ij}(x) =
\frac{\phi^{(\tilde{q})}_{ij}(x)}{\phi^{(\tilde{q})}_{ij}(x^{-1})},
\end{eqnarray*}

\noindent where the functions
$\phi^{(q)}_{ij}(x)$ and $\phi^{(\tilde{q})}_{ij}(x)$ are defined as follows,

\begin{eqnarray*}
\phi^{(q)}_{ij}(x) = \frac{\theta_q(xp^{A_{ij}/2})}{\theta_q(xq^{A_{ij}/2})},
~~~~~~\phi^{(\tilde{q})}_{ij}(x) =
\frac{\theta_{\tilde{q}}(xp^{A_{ij}/2})}{\theta_{\tilde{q}}(x\tilde{q}^{A_{ij}/2})}.
\end{eqnarray*}

\noindent These equations also imply that

\begin{eqnarray*}
\psi^{(q)}_{ij}(x)\psi^{(q)}_{ij}(x^{-1})=1,~~~~
\psi^{(\tilde{q})}_{ij}(x)\psi^{(\tilde{q})}_{ij}(x^{-1})=1.
\end{eqnarray*}

Eqs.(\ref{17}),(\ref{18}),(\ref{22}),(\ref{23})
give a level 1 representation of the algebra (\ref{39})-(\ref{51}) on the Fock
space of the Heisenberg algebra ${\cal H}_{q,p}(g)$.

To end this note, we would like to present some concluding remarks.
\begin{itemize}
\item With slight modifications to the zero mode contributions in the screening
currents of the quantum deformed $W$-algebras, we come to the new elliptic
deformed Kac-Moody algebra described in (\ref{39})-(\ref{51}).
Such algebras should be regarded
to be generated by the Fourier coefficients of the generating currents and
are associative algebras with unit. However, unlike the usual quantum deformations
of classical and affine Lie algebras, no Hopf algebraic or quasi-Hopf algebraic
structures can be defined over this new kind of elliptic algebras.
\item The algebra obtained in this note carries two deformation parameters $p, q$
in contrast to the quantum group and Yangian algebras where only one
deformation parameter $q$ ($\hbar$) is present. Moreover, the relations
for the ``positive'' and ``negative'' generating currents are deformed differently
in the sense that the deformation parameters are different in these relations.
The only known example of algebras of this kind before this note is the
algebras ${\cal A}_{\hbar,\eta}(\hat{g})$, which are members of the so-called
Hopf family of algebras \cite{g,f}. However, we are not able to
define the structure
of Hopf family over the present elliptic algebra. The relations between the
new elliptic algebra and the Hopf family of algebras will be an interesting
subject of further study.
\item For quantum groups and Yangian algebras, there exist different realizations
including the current realization, Reshetikhin-Semenov-Tian-Shansky realization
\cite{k} and the Drinfeld realizations. The elliptic algebra obtained in this note is
only realized in the current realization. It remains a hard problem to find the
other possible realizations, especially the realization which may have direct
relation with the quantum Yang-Baxter relations.
\item As mentioned earlier, the modification in the screening currents spoils
the feature that they commute with the quantum deformed $W$-algebra up to
total differences. Therefore, the relationship between the resulting algebra
and the quantum deformed $W$-algebra remains unclear. Presumably the quantum
deformed $W$-algebras can be obtained from our elliptic algebra in terms of
a $(q,p)$-deformed version of Hamiltonian reduction. It seems that this problem
deserves to be studied in detail. Notice that the $q$-difference version
of Hamiltonian reduction has recently been carried out by E.Frenkel, N.Reshetkhin,
M.A.Semenov-Tian-Shansky and A.V.Sevostyanov in \cite{aa},\cite{bb}.
\item The elliptic algebra is obtained here only with the level
$c=1$ bosonic representation. It seems that much more effort should be paid
towards the representation theory of this new kind of algebra, especially the
representations with higher level and the irreducibility, the spinor representations
and the level zero representations should be studied in detail. We do not
know whether there exist any relationship between the level zero form of our algebra
and the elliptic quantum groups proposed by G.Felder \cite{l}.
Hopefully there is some, but
this problem can only be answered after detailed study.
\item Last, the possible physical application of the new elliptic algebra should
be studied, e.g. whether there exist any physical models bearing the new algebra as
a symmetry. It is also interesting to mention that in conventional conformal field
theories, the screening currents of the usual $W$-algebras also do not close
into a single algebra. After similar modifications it is hopeful that they also
form closed algebras.
\end{itemize}

\end{document}